\documentclass[11pt,a4paper]{article}
\usepackage{soul}
\usepackage{amsmath}
\usepackage{amsfonts}
\usepackage[english]{babel}
\usepackage{fontenc}
\usepackage{euscript}
\usepackage{array}
\usepackage{vmargin}
\usepackage{amssymb}
\usepackage{latexsym}
\usepackage{enumerate}
\usepackage{graphicx,color}
\usepackage{hyperref}
\usepackage{commath}
\usepackage[titletoc,toc,title]{appendix}
\definecolor{c30}{rgb}{1,0.4,0.2}
\definecolor{c20}{rgb}{0, 0.1, 0.0}
\definecolor{c10}{rgb}{0.00,0,1.00}
\newcommand{\COM}[1]{}

\usepackage[utf8]{inputenc}

\def\rain{\to +\infty}

\def\convn{\xrightarrow [n\rain]{}}
\def\N{{\rm I\kern-.20em N}}
\def\R{{\rm I\kern-.20em R}}
\def\indi{{1\kern-.20em\rm I}}

\def\bkR{{\rm I\kern-.17em R}}
\def\bkN{{\rm I\kern-.20em N}}

\newtheorem{lemma}{Lemma}[section]
\newtheorem{theorem}{Theorem}[section]
\newtheorem{remark}{Remark}

\newtheorem{example}{Example}[section]

\newcommand{\bdem} {\begin{proof}}
\newcommand {\edem}{\hfill $\square$ \end{proof}}

\newcommand*{\affaddr}[1]{#1} % No op here. Customize it for different styles.
\newcommand*{\affmark}[1][*]{\textsuperscript{#1}}
\newcommand*{\email}[1]{\textit{#1}}

\newcommand{\PreserveBackslash}[1]{\let\temp=\\#1\let\\=\temp}
\newcolumntype{C}[1]{>{\PreserveBackslash\centering}p{#1}}
\newcolumntype{R}[1]{>{\PreserveBackslash\raggedleft}p{#1}}
\newcolumntype{L}[1]{>{\PreserveBackslash\raggedright}p{#1}}

\begin{document}
\title{Extremal behaviour of a periodically controlled sequence with imputed values}
\author{Helena Ferreira\affmark[1], Ana Paula Martins\affmark[1], Maria da Gra\c{c}a Temido\affmark[2] \\
\affaddr{\affmark[1]{\small Departamento de Matem\'atica, Centro de Matem\'atica e Aplica\c{c}\~oes (CMA-UBI),}}\\
\affaddr{\small Universidade da Beira Interior, Portugal}\\
\affaddr{\affmark[2] {\small CMUC, Departamento de Matem\'atica, Universidade de Coimbra, Portugal}}\\
\email{helenaf@ubi.pt, amartins@ubi.pt, mgtm@mat.uc.pt}
}
\date{}

  \maketitle

\begin{quote}
{\bf Abstract:}  Extreme events are a major concern in statistical modeling. Ran\-dom missing data can constitute a problem  when modeling such rare events. Imputation is crucial in these situations and therefore models that describe different imputation functions enhance possible applications and enlarge the few known families of models which cover these situations.  In this paper we consider a family of models $\{Y_n\},$ $n\geq 1,$ which can be associated to automatic systems which have a periodic control, in the sense that it is guaranteed that at instants multiple of $T,$ $T\geq 2,$ no value is lost. Random missing values are here replaced by the biggest of the previous observations up to the one surely registered. We prove that when the underlying sequence is stationary, $\{Y_n\}$ is $T-$periodic and if it also verifies some local dependence conditions then $\{Y_n\}$ verifies one of the well known $D^{(s)}_T(u_n),$ $s\geq 1,$ dependence conditions for $T-$periodic sequences. We also obtain the extremal index of $\{Y_n\}$ and relate it to the extremal index of the underlying sequence. A consistent estimator for the parameter that ``controls'' the missing values is here proposed and its finite sample properties are analysed.
The obtained results are illustrated with Markovian sequences of recognized interest in applications.

{\bf Resum\'e:}
Les événements extrêmes sont une préoccupation majeure dans la modélisation statistique. Les données absentes de façon aléatoire peuvent constituer un problème lors de la modélisation de ces événements rares. Les techniques d'imputation sont cruciaux dans ces situations et, par conséquent, les modèles qui décrivent différentes fonctions d'imputation améliorent les applications possibles et élargissent les quelques familles connues de modèles qui couvrent ces situations. Dans cet article, nous considérons une famille de modèles $ \{Y_n \}, $ $ n \geq 1, $ qui peuvent être associés à des systèmes automatiques qui ont un contrôle périodique, en ce sens qu'il est garanti qu'à des instants multiples de $ T, $ $ T \geq 2, $ aucune valeur n'est perdue. Les valeurs absentes de façon aléatoire sont ici remplacées par la plus grande des observations précédentes jusqu'à celle sûrement enregistrée. Nous prouvons que lorsque la séquence sous-jacente est stationnaire, $ \{Y_n \} $ est $ T- $ périodique et si elle vérifie encore certaines conditions de dépendance locales, alors $ \{Y_n \} $ vérifie l'une des conditions de dépendence $ D^{ (s)}_T (u_n), $ $ s\geq 1, $ pour des séquences $ T- $ périodiques. Nous obtenons également l'indice extrémal de $ \{Y_n \} $ que nous relions à l'indice extrémal de la séquence sous-jacente. Un estimateur convergent pour le paramètre qui ``contrôle'' les valeurs absentes est ici proposé et ses propriétés à distance fini sont analysées. Les résultats obtenus sont illustrés en utilisant des séquences markoviennes qui ont un intérêt reconnu dans les applications.

\end{quote}\vspace{0.3cm}

%% --------------------------------------------------------------
%% Inicio do texto do artigo
%% --------------------------------------------------------------

\section{Introduction and preliminary results}
\label{intro}

Data collection is prevalent in everyday life and is used in several do\-mains, such as finance, climate observation, computer science, etc. The main goal of any data collection effort is to compile quality data, but issues with missing data oftenly occur  when data is measured and recorded. The data unavailability may be caused by the failure of some system, such  as a reading device, or simply by lack of retention due to the intrinsic properties of the data, {\it{e.g.}} financial or environmental data only reported at certain time instants (Hall and H\"usler (2006); Hall and Scotto (2008) and references therein).

Many analysis methods require for the use of imputation, {\it{i.e.}} missing values to be replaced with reasonable values up-front. An overview about univariate time series imputation can be found in Moritz {\it{et al.}} (2015)  and an introduction to \textsf{R}'s package \textsf{imputeTS}, which is solely dedicated to univariate time series imputation is presented in  Moritz and Bartz-Beielstein (2017).

Falk {\it{et al.}} (2011) summarize the several strategies that are usually applied when missing data values occur in time series: (i) the missing value is replaced by a predefined value $x_0$ which can be sometimes 999 (if one is interested in small values and no such large values occur) or -1 (if one is  interested in  large values and no negative values occur); (ii) the data is completely lost and the time series is sub-sampled with a smaller (and random) sample size; (iii) an automatic measurement device is used to replace  the missing data by a proxy value.

The extremal properties of sequences with random missing values replaced by 0 were studied by Falk {\it{et al.}} (2011). The sub-sample referred in strategy (ii) above, may result from missing values that occur according to some deterministic pattern or occur randomly. The effect of deterministic missing values on the properties of strictly stationary (stationary) sequences has been studied by  Ferreira and Martins (2003), Martins and Ferreira (2004), Scotto {\it{et al.}} (2003), among others. Random missing values have been considered by Weissman and Cohen (1995) for the case of constant failure probability and independent failures, and their results were generalized for situations where the failure pattern has a weak dependence structure by  Hall and H\"usler (2006). This was pursued by Hall and Scotto (2008), when the underlying process is represented as a moving average driven by heavy-tailed innovations and the sub-sampling process is strongly mixing.

When the missing values are replaced by an automatic measurement device the resulting sample will be a mixture of two original samples. This case was considered by   ulting sample will be a mixture of two original samples. This case was considered by  Hall and H\"usler (2006) and later by  Hall and Temido (2009) in the context of max-semistability. There they discussed the extremal properties of a model where missing values are replaced by independent replicas of the original values. Investigating the extremal properties of models that describe other imputation functions enhances possible applications in situations where it may be of interest to avoid the occurrence of missing values and an automatic replacement of a device or machine may be available. This situation motivated the model we consider in this paper and that we describe in what follows.\vspace{0.3cm}

Let us consider a system with a periodic control, in the sense that it is guaranteed that at instants multiple of $T,$ $T\geq 2,$ no value is lost. If for some reason there are missing values, it is then natural to use the observations that are for surely registered, due to the periodic control, in the replacement of these values.  A model that translates this idea, where a missing value is replaced by the biggest of the previous observations up to the one at an instant multiple of $T,$ which is surely registered, is defined by
\begin{equation}\label{1.1}
Y_n=U_nX_n+(1-U_n)\bigvee_{i=n-1}^{\left[\frac{n-1}{T}\right]T} U_iX_i\quad n \geq 1,\ T\geq 2,
\end{equation}
 where  $[a]$ denotes the integer part of $a\in \R,$  $\{U_n\}_{n\geq 0}$ is a sequence of independent variables, such that, for all $k\geq 0,$ $U_{kT}=1$ almost surely, and  $U_n$ follows a Bernoulli distribution with parameter $p\in ]0,1[,$ for all $n \neq kT.$  $\{X_n\}_{n\geq 0}$ denotes a positive stationary sequence, independent of  $U_n,$ $n \geq 0,$ with marginal continuous distribution function (d.f.) $F.$ Here and throughout  $\bigvee_{j=s}^t W_j\equiv \bigvee_{j=t}^s W_j,$ for any $s,t \in \mathbb{Z},$ $s<t,$ and random variables $W_j.$

Model (\ref{1.1}) can be associated to automatic systems which have a manual periodic verification. As we can see, at each instant $n\neq kT,\, k \geq 0,$ we can observe $X_n$ or in the case that it is not observed it is replaced by the maximum of the previous observations up to the last one that was surely registered. The registration of the observations is periodically controlled, with the guarantee that at instants $kT,$ $k\geq 0,$ no observations are lost and therefore in the period $\left\{\left[\frac{n-1}{T}\right]T,\ldots, n-1\right\}$ at least the observation with index $\left[\frac{n-1}{T}\right]T$  is available. The case $T=1$, if considered, would correspond to the non occurrence of missing data.

To better understand model (\ref{1.1}) let us consider the following illustrative example.

\begin{example}\label{ex1}
{\rm{Let  $\{Z_n\}_{n\geq -1}$ be a sequence of independent and identically distributed (i.i.d.) random variables with unit Fr\'echet marginal d.f. $F_Z$. With this sequence we define the moving maxima  $X_n=\frac{1}{2}Z_n\vee \frac{1}{2}Z_{n-1},\ n\geq 0,$ which is stationary and also has unit Fr\'echet margins.

A model with the characteristics of (\ref{1.1}) is given by
\begin{equation}\label{1.2}
Y_n=U_nX_n+(1-U_n)U_{n-1}X_{n-1},\quad n\geq 1,
\end{equation}
where $\{U_n\}_{n\geq 0}$ is an independent Bernoulli sequence with parameter $p$ for all $n\neq 2k$ and $U_{2k}=1,$ $k\geq 0,$ almost surely.

 Here  $\{Y_n\}$ is controlled at instants which are multiples of 2 ($T=2$), so at these instants we always have the guarantee that an observation of the moving maxima was retained. At all the  other instants we can have a moving maxima value, whenever it is observed, or the maximum of the previous observations up to the last observation ``controlled'', which in this case corresponds only to the previous observation because $T=2.$ In Figure \ref{fig1} this becomes clear with 100 observation of (\ref{1.2}), since newly imputed observations are marked differently than the rest of the series and the instants $2s,$ $s\geq 1,$ are highlighted.

 \begin{figure}[!ht]\label{fig1}
 \begin{center}
 \includegraphics[scale=0.35]{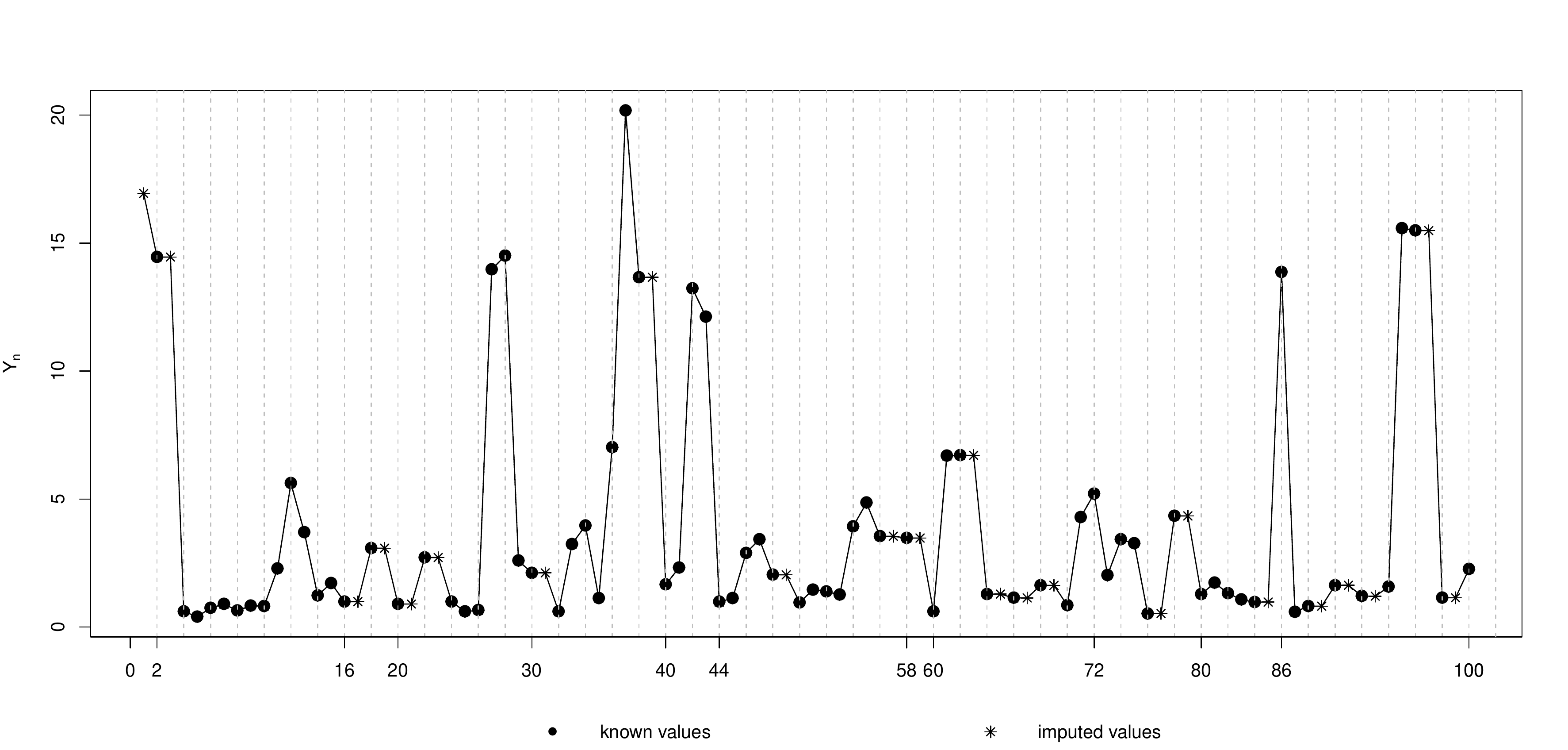}
 \caption{100 observations of model  (\ref{1.2}) with $p=0.5$}
  \end{center}
 \end{figure}

 We shall return to this simple example, throughout the work, to illustrate several of the results presented.}}
\end{example} \vspace{0.3cm}

Our main goal is to characterize the extremal behaviour of  $\{Y_n\}_{n\geq 1}$ given in (\ref{1.1}). In order to achieve this, in Section 3, lets start by noting that, for all $ k \geq 0$ and for all $x\in \R,$ the  marginal d.f.'s of $\{Y_n\}_{n\geq 1}$ satisfy
 $$ F_{Y_{kT+j}}(x)=
\begin{cases}  F(x), & \,\, j \in \{0,1\}, T \geq 2,  \\ pF(x)+(1-p)G_j(x), & \,\,  j \in \{2,\ldots,T-1\},\,\,T \geq 3, \end{cases}$$
 with $kT+j \neq 0$ and
 \begin{eqnarray}\label{1.3}
\lefteqn{\hspace{-1cm}G_j(x)=\sum_{\varnothing \subseteq S\subseteq \{kT+1,\ldots,kT+j-1\}}p^{|S|}(1-p)^{j-1-|S|}\times}\\
&\times& F_{kT, kT+1,\ldots,kT+j-1}
(x,x(\delta_{kT+1}(S))^{-1},\ldots,x(\delta_{kT+j-1}(S))^{-1}),\hspace{-2cm}\nonumber
 \end{eqnarray}
 where $\delta_{\bullet}(B)$ denotes the Dirac measure on $B,$ $F_{i_1,\ldots,i_p}$ denotes the d.f. of the vector $(X_{i_1},\ldots,X_{i_p})$ and we conventione that $x\times \frac{1}{0}=+\infty.$

We also prove that the sequence $\{Y_n\}_{n\geq 1}$ defined in (\ref{1.1}) is a $T$ periodic sequence. Indeed, for any choice of integers $1\leq i_1<i_2<\ldots <i_p,$ we have\vspace{-0.1cm}
 \begin{eqnarray*}
 \lefteqn{P\left(\bigcap_{s=1}^p \{Y_{i_s+T}\leq x_s\}\right)}\\[0.1cm]
 &=& P\left(\bigcap_{s=1}^p \left\{U_{i_s+T}X_{i_s+T}+(1-U_{i_s+T})\bigvee_{j=i_s+T-1}^{\left[\frac{i_s+T-1}{T}\right]T}U_jX_j\leq x_s\right\}\right)\\[0.1cm]
 &=&E\left(P\left(\bigcap_{s=1}^p \left\{U_{i_s+T}X_{i_s+T}+(1-U_{i_s+T})\bigvee_{j=i_s+T-1}^{\left[\frac{i_s+T-1}{T}\right]T}U_jX_j\leq x_s\right\}\right)\Biggm|\right.\\&&\hspace{7cm}\left. U_{\left[\frac{i_s+T-1}{T}\right]T},\ldots,U_{i_s+T}\right)\\[0.1cm]
 &=&E\left(P\left(\bigcap_{s=1}^p \left\{U_{i_s}X_{i_s+T}+(1-U_{i_s})\bigvee_{j=i_s+T-1}^{\left[\frac{i_s+T-1}{T}\right]T}U_jX_j\leq x_s\right\}\right)\Biggm| \right.\\&&\hspace{8cm}\left. U_{\left[\frac{i_s-1}{T}\right]T},\ldots,U_{i_s}\right)\\[0.2cm]
 &=&E\left(P\left(\bigcap_{s=1}^p \left\{U_{i_s}X_{i_s}+(1-U_{i_s})\bigvee_{j=i_s-1}^{\left[\frac{i_s-1}{T}\right]T}U_jX_j\leq x_s\right\}\right)\Biggm|\right.\\&&\hspace{8cm}\left.  U_{\left[\frac{i_s-1}{T}\right]T},\ldots,U_{i_s}\right),
 \end{eqnarray*}
 since $\{U_n\}_{n\geq 1}$ is a $T-$periodic sequence, $\{X_n\}_{n\geq 0}$ is a stationary sequence and they are independent. From now on we use the notation $F_j(x):=F_{Y_{j}}(x)$, for any $j\geq 1$, for which holds $F_j(x):=F_{j+kT}(x),$ for all $k\geq 1$ and $T\geq 2$.\vspace{0.3cm}

Extreme value theory known for periodic sequences can then be
applied to this periodically controlled sequence with imputed values $\{Y_n\}_{n\geq
1},$ since it is also  a $T-$periodic sequence. Alpuim (1988)
showed that under Leadbetter's global mixing condition $D(u_n),$ the
only possible limit laws for the normalized maxima of a
$T-$periodic sequence are the three extreme value distributions and generalized, as well, the definition of extremal index for such sequences.
Under local mixing conditions $D_T^{(s)}(u_n),\ s=1,2,$  Ferreira (1994)
studied the extremal behaviour of periodic sequences and under
the weaker local mixing conditions $D_T^{(s)},\ s\geq 3,$  Ferreira
and Martins (2003) obtained the expression for the extremal index
of a $T$-periodic sequence from the joint distribution of $s$
consecutive variables of the sequence.

In Section 3 we obtain necessary conditions, that rely on the underlying sequence $\{X_n\}_{n\geq 0},$ for sequence $\{Y_n\}_{n\geq 1}$ to satisfy Leadbetter's $D(u_n)$ condition, as well as  some local dependence condition $D_T^{(s)},$ $s\geq 1,$ for $T-$periodic sequences. The validation of these conditions will permit the determination of its extremal index expression. The results here obtained are illustrated with examples of recognized interest in applications, such as Markovian sequences.

The next section is devoted to the estimation of the model parameter $p\in ]0,1[.$ We propose a consistent estimator for this parameter and analyse its finite sample behavior.

\section{Model parameter estimation}
\label{sec:1}

The proposed model (\ref{1.1}) depends on an unknown parameter $p\in ]0,1[,$ that ``controls'' the number of missing values, and on an underlying  stationary sequence with unknown marginal d.f. $F.$ The estimation of $p$ and $F$ is therefore essential for practical applications of this model.

The next result, that characterizes the probabilities $P\left(\bigcap_{j=1}^{T-1}\{Y_{sT+j}=Y_{sT}\}\right),$ $s\geq 1,$ $T\geq 2,$   gives a simple procedure to estimate the parameter $p$ involved in the definition of model (\ref{1.1}).

\begin{theorem} \label{teo1}
For the sequence $\{Y_n\}_{n\geq 1}$ defined in (\ref{1.1}) it holds
$$P\left(\bigcap_{j=1}^{T-1}\{Y_{sT+j}=Y_{sT}\}\right)=(1-p)^{T-1},\quad s\geq 1, T\geq 2.$$
\end{theorem}
\bdem
For all $s\geq 1,$  $Y_{sT}=X_{sT}$ almost surely and, if $U_{sT+j}=1$ for some  $j\in \{1,\ldots,T-1\}$, we have $P\left(\bigcap_{j=1}^{T-1}\{Y_{sT+j}=Y_{sT}\}\right)=0,$ because the underlying d.f.'s are continuous. Therefore we can write\vspace{-0.1cm}
\begin{eqnarray}\label{2.1}
\lefteqn{P\left(\bigcap_{j=1}^{T-1}\{Y_{sT+j}=Y_{sT}\}\right)}\nonumber\\&=& P\left(\bigcap_{j=1}^{T-1}\left\{U_{sT+j}X_{sT+j}+(1-U_{sT+j})\bigvee_{i=sT+j-1}^{\left[\frac{sT+j-1}{T}\right]T}U_iX_i=X_{sT}\right\}\right)\nonumber \\
&=&(1-p)^{T-1}P\left(\bigcap_{j=1}^{T-1}\left\{\bigvee_{i=sT+j-1}^{\left[\frac{sT+j-1}{T}\right]T}U_iX_i=X_{sT}\right\}\ \Biggl|\ \bigcap_{j=1}^{T-1}U_{sT+j}=0\right).
\end{eqnarray}
Now since $U_{sT+j}=0,$ for all $j=1,\ldots,T-1,$ the several maxima in (\ref{2.1}) are all equal to the variable $X_{sT}$ and the result follows immediately.
\edem\vspace{0.2cm}

The way to estimate  parameter $p$ of model (\ref{1.1}) becomes clear from the previous result. So, if $(Y_1,\ldots,Y_n)$ is a random sample of $\{Y_n\}_{n\geq 1}$  an estimator for $p\in ]0,1[$ is given by\vspace{-0.2cm}
\begin{equation}\label{2.2}
\widehat{p}_{n,T}=1-\left(\frac{1}{\left[\frac{n+1}{T}\right]-1}\sum_{s=1}^{\left
[\frac{n+1}{T}\right]-1} \indi_{A_s}\right)^{\frac{1}{T-1}},\ T\geq 2,
\end{equation}\vspace{-0.2cm}
with $A_s=\displaystyle{
\bigcap_{j=1}^{T-1}\{Y_{sT+j}=Y_{sT}\}}.$

From the weak law of large numbers and the fact that \linebreak $E[\indi_{\{\bigcap_{j=1}^{T-1}\{Y_{sT+j}=Y_{sT}\}\}}]=P\left(\bigcap_{j=1}^{T-1}\{Y_{sT+j}=Y_{sT}\}\right),$ we can state that estimator (\ref{2.2}) is a consistent estimator for $p.$

The d.f. $F$ can be estimated from the observations $Y_{sT}\equiv X_{sT},$ $s\geq 1,$ with the empirical d.f. or a kernel estimator. A review of these estimators and an explanation on their functionality and applicability in \textsf{R} can be found in  Quintela-del-R\'{\i}o and Est\'evez-P\'erez (2012).

\subsection{Simulation results}
\label{sec:2}

We now analyze the finite sample properties of the estimator given in (\ref{2.2}) with simulated data from model (\ref{1.2}) given in Example \ref{ex1}. Each simulated data set consists of 1000 independent copies of $n$ realizations of a random sample $(Y_1,\ldots,Y_n)$ of (\ref{1.2}) having one particular value of $p\in]0,1[$ out of five, for   $T=2.$ Three different sample sizes are considered for each data set. The sample means $\widehat{\mu}(\widehat{p}_{n,2})$ and the sample standard deviations $\widehat{\sigma}(\widehat{p}_{n,2})$ of the estimates $\widehat{p}_{n,2,i},$ $i=1,\ldots,1000,$ depending on the sample size $n,$ were computed. The bias and the root mean squared errors ($RMSE(\widehat{p}_{n,2})$) were also determined. Table \ref{tab1} summarizes the estimation results obtained. The estimator has a good behavior even for small sample sizes.

\begin{table}[t!]\label{tab1}
\caption{Various statistical results for the estimation of $p\in]0,1[$ in model (\ref{1.2})}
\begin{center}
\begin{tabular}{lL{2cm}crC{2cm}c}
\hline
 $p$ & n& $\widehat{\mu}(\widehat{p}_{n,2})$ & $BIAS(\widehat{p}_{n,2})$ & $\widehat{\sigma}(\widehat{p}_{n,2})$&$RMSE(\widehat{p}_{n,2})$\\ \hline
 0.10& $n=250$&0.0989&-0.0011&  0.0267&  0.0267\\
 & $n=1000$&0.1008&0.0008& 0.0137& 0.0137\\
 & $n=5000$&0.1001& 0.0001&0.0059&0.0059\\
 0.25& $n=250$& 0.2492&-0.0008&0.0395&0.0395\\
 & $n=1000$&0.2494&-0.0006&0.0198&0.0198\\
 & $n=5000$&0.2495&-0.0005&0.0083&0.0083\\
 0.50& $n=250$& 0.5025&0.0025&0.0461&0.0461\\
 & $n=1000$&  0.4997&-0.0003&0.0222&0.0223\\
 & $n=5000$&  0.4999& 0.0001&0.0096&0.0096\\
 0.75& $n=250$& 0.7510&0.0010&0.0373&0.0373\\
 & $n=1000$&  0.7506&0.0006&0.0188&0.0188\\
 & $n=5000$&  0.7498& -0.0002&0.0085&0.0085\\
 0.90& $n=250$& 0.9000&0.0000&0.0263&0.0263\\
 & $n=1000$&   0.9001&0.0001& 0.0129& 0.0129\\
 & $n=5000$&  0.9001&0.0001&0.0060&0.0060\\
  \hline
 \end{tabular}
\end{center}
\end{table}

\section{Extremal behaviour}

The study of the extremal behavior of stationary or periodic sequences most often relies upon the verification of appropriate dependence conditions which assure that the limiting distribution of the ma\-xi\-mum term is of the same type as the limiting distribution of the maximum of i.i.d. random variables. Usual conditions used in the literature for stationary and periodic sequences are Leadbetter's $D(u_n)$ condition (Leadbetter 1974), conditions $D^{(s)}(u_n),$ $s\geq 1,$ of Chernick {\it{et al.}} (1991) and conditions $D^{(s)}_T(u_n),$ $s, T\geq 1,$ of Ferreira (1994) and  Ferreira and Martins (2003).\vspace{0.3cm}

For any sequence of real numbers $\{u_n\}_{n\geq 1},$ the dependence condition $D(u_n),$ for the sequence $\{Y_n\}_{n\geq 1},$ of Leadbetter,
states that $\alpha_{n,\ell_n} \to 0,$ as $n\to +\infty$, for some sequence $\ell_n = o(n)$, where
\begin{eqnarray*}
\alpha_{n,\ell_n} &=&\sup\{|P(M^{(\textbf{Y})}_{i_1,i_1+p}\leq  u_n, M^{(\textbf{Y})}_{j_1,j_1+q} \leq u_n) -P(M^{(\textbf{Y})}_{i_1,i_1+p} \leq u_n)\times\\ &&\times P(M^{(\textbf{Y})}_{j_1,j_1+q}\leq u_n)| :1 \leq i_1 < i_1 + p + \ell_n \leq j_1 < j_1 + q \leq n\},
\end{eqnarray*}
and $M^{(\textbf{Y})}_{s,s+t}$ denotes $\displaystyle{\bigvee_{i=s}^{s+t}Y_i} \equiv \displaystyle{\bigvee_{i=s+t}^{s}Y_i}.$ \vspace{0.3cm}

In the next result we show that if condition  $D(u_n)$ holds for the stationary sequence $\{X_n\}_{n\geq 0}$ then it also holds for the $T-$ periodic sequence $\{Y_n\}_{n\geq 1}.$

\begin{theorem}\label{teo2}
If, for any positive sequence  $\{u_n\}_{n\geq 1},$ condition $D(u_n)$ holds for the stationary sequence $\{X_n\}_{n\geq 0}$ then it also holds for the sequence $\{Y_n\}_{n\geq 1}$ defined in (\ref{1.1}).
 \end{theorem}
\bdem
Consider the sets of consecutive integers $A_p= \left\lbrace i_1, \cdots, i_1 +p \right\rbrace$ and  $B_q= \left\lbrace j_1, \cdots, j_1 + q \right\rbrace$, with $j_1-(i_1+p) > \ell_n.$
	
The definition of $\{Y_n\}_{n\geq 1}$ induces over  $\{X_n\}_{n \geq 0}$ the corresponding sets of integers  \linebreak $A^{'}_p= \left\lbrace \left[\frac{i_1-1}{T}\right] T, \cdots, i_1 +p \right\rbrace$ and $B^{'}_q=\left\lbrace \left[\frac{j_1-1}{T}\right] T, \cdots, j_1 + q \right\rbrace$, such that  $\left[\frac{j_1-1}{T}\right] T - (i_p+p) > \ell_n-T.$
Hence, each realization $u$ of $\textbf{U}:=\left\{U_i, \ i \in A_p^{'} \cup B_q^{'} \right\}$  gives rise to another pair of subsets of positive integers, say $A_p^{(\textbf{U})}$ and $B_q^{(\textbf{U})}$, and therefore
\begin{eqnarray*}
\lefteqn{\hspace{-1cm}P\left(M^{(\textbf{Y})}(A_p)\leq u_n, \ M^{(\textbf{Y})}(B_q)\leq u_n\right)}\nonumber\\[0.3cm]
	&=& E\left(P\left(M^{(\textbf{X})}(A_p^{(\textbf{U})})\leq u_n, \ M^{(\textbf{X})}(B_q^{(\textbf{U})})\leq u_n \mid \textbf{U}\right) \right).
\end{eqnarray*}

	Due to the fact that $\{X_n \}_{n\geq 0}$ satisfies $D(u_n)$ condition, the last average becomes
	\begin{equation*}
	\begin{aligned}
\lefteqn{E\left(P\left(M^{(\textbf{X})}(A_p^{(\textbf{U})})\leq u_n\mid \ \textbf{U}\right)P\left(M^{(\textbf{X})}(B_q^{(\textbf{U})})\leq u_n \mid \textbf{U}\right)\right) + \mathcal{O}(\alpha_{n, \ell_n-T})} \\[0.3cm]
	&=E\left(P\left(M^{(\textbf{X})}(A_p^{(\textbf{U})})\leq u_n \mid \textbf{U} \right)\right)E\left(P\left(M^{(\textbf{X})}(B_q^{(\textbf{U})})\leq u_n \mid  \textbf{U}\right) \right) + \mathcal{O}(\alpha_{n, \ell_n-T})
	\end{aligned}
	\end{equation*}
because $\{U_n\}_{n\geq 0}$ is a sequence of independent variables.

Returning to the sequence $\{Y_n\}_{n\geq 1}$ we deduce
	\begin{eqnarray*}
\lefteqn{|P\left(M^{(Y)}(A_p)\leq u_n, \ M^{(Y)}(B_q)\leq u_n\right) -}\\
 &&\qquad \qquad -P\left(M^{(Y)}(A_p)\leq u_n\right)P\left( M^{(Y)}(B_q)\leq u_n\right)|\leq \beta_{n, \ell_n},
\end{eqnarray*}
with $\beta_{n, \ell_n}=\alpha_{n, \ell_n-T}$ and for $\ell'_n=\ell_n+T=o(n)$ we have $\beta_{n,\ell'_n}\convn 0,$ as required.\edem\vspace{0.3cm}

If the $T$-periodic sequence $\{Y_n\}_{n\geq 1}$ satisfies condition $D(u_n),$ for all $\{u_n\}_{n\geq 1},$ then we say that it also satisfies  condition $D_T^{(s)}(u_n),$ $T, s\geq 1,$  when
there exists a sequence of integers $\{k_n\}_{n\geq 1}$
such that $k_n\to +\infty,$ $
k_n\frac{\ell_n}{n}\to 0,$ $k_n \alpha_{n,\ell_n}\to 0,$ as $n\to +\infty,$ and
\begin{equation}\label{3.1}
\lim_{n\to
+\infty}n\frac{1}{T}\sum_{i=1}^{T}P\left(Y_i>u_n\geq M^{(\bf{Y})}_{i+1,i+s-1},\ M^{(\bf{Y})}_{i+s,\left[\frac{n}{k_nT}\right]T}>u_n\right)=0,
\end{equation}

\noindent where
$M^{(\bf{Y})}_{i,j}=-\infty,$ for $i>j,$ and $M^{(\bf{Y})}_{i,i}=Y_i.$ These local dependence conditions were first defined in Ferreira (1994), for $s=1$ and $s=2.$ This family was later enlarged by  Ferreira and Martins (2003) with values of  $s\geq 3.$\vspace{0.3cm}

Observe that, when $s\geq 2,$ condition (\ref{3.1}) is implied
by%\vspace{0.3cm}
\begin{equation*}\label{3.2}
\lim_{n\to
+\infty}n\frac{1}{T}\sum_{i=1}^{T}\sum_{j=i+s}^{[\frac{n}{k_nT}]T}P(Y_i>u_n,Y_{j-1}\leq
u_n<Y_j)=0,
\end{equation*}
which limits the distance between
exceedances of level $u_n,$ {\it{i.e.}}, in each interval there can
only occur more than one exceedance of $u_n$ if separated by less
than $s-1$ non-exceedances of $u_n.$ Consequently, the local
dependence conditions $D^{(s)}_T,\ s\geq 1,$ become weaker as the
value of $s$ increases and thereby enhance the number of processes to which our results apply.  Condition (\ref{3.1}), when $T=1,$ coincides with the one considered in $D^{(s)}(u_n)$ of  Chernick {\it{et al.}} (1991) for stationary sequences.  \vspace{0.3cm}

Under some local dependence condition $D_T^{(s)}(u_n),$ $s\geq 1,$ the extremal
index of $\{Y_n\}_{n\geq 1}$, $\theta_{\bf{Y}},$ is defined by  $-\frac{1}{\tau}\displaystyle{\log(\lim_{n\to
\infty} P(M^{(\bf{Y})}_{1,n} \leq u_n))},$ where $\{u_n\}_{n\geq 1}$ is a sequence of positive numbers such that
\begin{equation*}
\tau=\lim_{n\to +\infty}n\frac{1}{T}\sum_{j=1}^T
\overline{F}_j(u_n)=\lim_{n\to +\infty}n\left(1-\frac{1}{T}\sum_{j=1}^T F_j(u_n)\right),
\end{equation*}
with $\overline{F}_j(x)=1-F_j(x)$  the tail functions of $Y_j$ for $j=1,...,T.$ Hence, it can be computed from\vspace{0.3cm}
\begin{equation}\label{3.3}
\theta_{{\bf{Y}}}=\frac{1}{\tau}\lim_{n\to
+\infty}\displaystyle{n\frac{1}{T}\sum_{i=1}^T P(Y_i>u_n\geq
M^{(\bf{Y})}_{i+1,i+s-1})}.
\end{equation}\vspace{0.3cm}

In order to apply the previous results we shall impose the following two conditions on the tail of the d.f.'s $F$ and $G_j,$ $j=2,\ldots,T-1, \, T\geq 3,$ where $G_j$ denotes the d.f. of $\bigvee_{i=j-1}^{\left[\frac{j-1}{T}\right]T}U_jX_j$, defined by (\ref{1.3}). Namely
\begin{equation}\label{3.4}
n\overline{F}(u_n)\convn \tau_{\bf{X}}\qquad \textrm{and}\qquad n\overline{G}_j(u_n)\convn \tau_j.
\end{equation}
The particular case $T=2$ leads to $\tau=\tau_{\bf{X}}$. Under such conditions, for $T \geq 3$, we have
\begin{eqnarray*}
\lefteqn{n\frac{1}{T}\sum_{j=1}^T\overline{F}_j(u_n)=} \nonumber\\[0.3cm]
&=&n\frac{1}{T}\left\{\sum_{j=2}^{T-1}
\left(p\overline{F}(u_n)+(1-p)P\left(\bigvee_{i=j-1}^{\left[\frac{j-1}{T}\right]T}U_i X_i>u_n\right)\right)+2\overline{F}(u_n)\right\}\nonumber\\[0.3cm]
\end{eqnarray*}
\begin{eqnarray}\label{3.5}
&=& n\frac{1}{T}((T-2)p+2)\overline{F}(u_n)+(1-p)\frac{n}{T}\sum_{j=2}^{T-1}P\left(\bigvee_{i=j-1}^{\left[\frac{j-1}{T}\right]T}U_iX_i>u_n\right)\nonumber\\[0.3cm]
&=&\frac{1}{T}((T-2)p+2)n\overline{F}(u_n)+(1-p)\frac{1}{T}\sum_{j=2}^{T-1}n\overline{G}_j(u_n)\nonumber\\[0.3cm]
&\convn& \frac{1}{T}\left\{((T-2)p+2)\tau_{\bf{X}}+(1-p)\sum_{j=2}^{T-1}\tau_j\right\}=\tau.
\end{eqnarray}\vspace{0.3cm}

We derive now a relation between a slightly stronger condition than $D^{(T+1)}(u_n)$ of Chernick {\it{et al.}} (1991) for the underlying stationary sequence $\{X_n\}_{n\geq 0}$ and condition $D^{(T+1)}_T(u_n)$ for the $T$-periodic sequence $\{Y_n\}_{n\geq 1}$.

\begin{theorem}\label{teo3}
If for any positive sequence $\{u_n\}_{n\geq 1},$ the stationary sequence $\{X_n\}_{n\geq 0}$ satisfies
\begin{equation} \label{3.6}
nP\left(\bigvee_{j=0}^T X_j>u_n\geq X_{T+1},\ M^{(\bf{X})}_{T+2,\left[\frac{n}{k_nT}\right]T+T}>u_n\right)\convn 0
\end{equation}
 then condition $D_T^{(T+1)}(u_n)$ holds for the sequence $\{Y_n\}_{n\geq 1}.$
\end{theorem}
\bdem
Observe first that
\begin{eqnarray*}
\lefteqn{\hspace{-1cm}n\frac{1}{T}\sum_{i=1}^{T}P\left(Y_i>u_n\geq M^{(\bf{Y})}_{i+1,i+T},\ M^{(\bf{Y})}_{i+T+1,\left[\frac{n}{k_nT}\right]T}>u_n\right)}\\[0.3cm]
& \leq & n\frac{1}{T}\sum_{i=1}^{T}P\left(Y_i>u_n\geq Y_{i+1},\ M^{(\bf{Y})}_{i+T+1,\left[\frac{n}{k_nT}\right]T}>u_n\right)\\[0.3cm]
&\leq & n\frac{1}{T}\sum_{i=1}^{T}P\left(Y_i>u_n\geq Y_{i+1},\ M^{(\bf{X})}_{i+2,\left[\frac{n}{k_nT}\right]T}>u_n\right).
\end{eqnarray*}
In what concerns the probability involved in this last sum, for\linebreak $i \,\in\, \{1,2,\ldots,T-2 \}$, we have
\begin{eqnarray*}
\lefteqn{\hspace{-1cm}P\left(Y_i>u_n\geq Y_{i+1},\ M^{(\bf{X})}_{i+2,\left[\frac{n}{k_nT}\right]T}>u_n\ \biggl| \  U_i=1,U_{i+1}=1\right)}\\[0.3cm]
& = &P\left(X_i>u_n\geq X_{i+1},\ M^{(\bf{X})}_{i+2,\left[\frac{n}{k_nT}\right]T}>u_n\right)\\[0.3cm]
& =& P\left(X_T>u_n\geq X_{T+1},\ M^{(\bf{X})}_{T+2,\left[\frac{n}{k_nT}\right]T+T-i}>u_n\right)\\[0.3cm]
&\leq& P\left(\bigvee_{j=0}^T X_j>u_n\geq X_{T+1},\ M^{(\bf{X})}_{T+2,\left[\frac{n}{k_nT}\right]T+T}>u_n\right),
\end{eqnarray*}
as well as
\begin{eqnarray}\label{3.7}
\lefteqn{\hspace{-1cm}P\left(Y_i>u_n\geq Y_{i+1},\ M^{(\bf{X})}_{i+2,\left[\frac{n}{k_nT}\right]T}>u_n\ \biggl| \  U_i=0,U_{i+1}=1\right)}\nonumber\\[0.3cm]
 &=& P\left(\bigvee_{j=0}^{i-1} U_jX_j>u_n\geq X_{i+1},\ M^{(\bf{X})}_{i+2,\left[\frac{n}{k_nT}\right]T}>u_n\right)\nonumber\\[0.3cm]
&\leq& P\left(\bigvee_{j=0}^T X_j>u_n\geq X_{T+1},\ M^{(\bf{X})}_{T+2,\left[\frac{n}{k_nT}\right]T+T}>u_n\right),
\end{eqnarray}
where we used the stationarity of $\{X_n\}$. The remaining probabilities are equal to zero, since
\begin{eqnarray*}
\lefteqn{\hspace{-1cm}P\left(Y_i>u_n\geq Y_{i+1},\ M^{(\bf{X})}_{i+2,\left[\frac{n}{k_nT}\right]T}>u_n\ \biggl| \ U_i=1,U_{i+1}=0\right)}\\[0.3cm]
&& = P\left(X_i>u_n\geq \bigvee_{j=0}^{i} U_j X_j,\ M^{(\bf{X})}_{i+2,\left[\frac{n}{k_nT}\right]T}>u_n \biggl| \ U_i=1\right)\\[0.3cm]
&&\leq P\left(X_i>u_n \geq  X_i\right)=0
\end{eqnarray*}
and
\begin{eqnarray*}
\lefteqn{P\left(Y_i>u_n \geq Y_{i+1},\ M^{(\bf{X})}_{i+2,\left[\frac{n}{k_nT}\right]T}>u_n\ \biggl| \ U_i=0,U_{i+1}=0\right)}\\[0.3cm]
&& = P\left( \bigvee_{j=0}^{i-1} U_j X_j>u_n \geq \bigvee_{j=0}^{i} U_j X_j,\ M^{(\bf{X})}_{i+2,\left[\frac{n}{k_nT}\right]T}>u_n | U_i=0\right)\\[0.3cm]
&&\leq P\left( \bigvee_{j=0}^{i-1} U_j X_j>u_n \geq \bigvee_{j=0}^{i-1} U_j X_j \right)=0.
\end{eqnarray*}
For $i=T-1$ and $i=T$ we respectively have
\begin{eqnarray*}
\lefteqn{\hspace{-1cm}P\left(Y_{T-1}>u_n \geq Y_{T},\ M^{(\bf{X})}_{T+1,\left[\frac{n}{k_nT}\right]T}>u_n\right)}\\[0.3cm]
& \leq&  P\left(X_{T-1}>u_n \geq X_{T},\ M^{(\bf{X})}_{T+1,\left[\frac{n}{k_nT}\right]T}>u_n\right)\\[0.3cm]
&& + P\left( \bigvee_{j=0}^{T-2}  X_j>u_n \geq  X_T,\ M^{(\bf{X})}_{T+1,\left[\frac{n}{k_nT}\right]T}>u_n \right)\\[0.3cm]
&\leq&  2 P\left( \bigvee_{j=0}^{T}  X_j>u_n \geq  X_{T+1},\ M^{(\bf{X})}_{T+2,\left[\frac{n}{k_nT}\right]T+T}>u_n \right)
\end{eqnarray*}
and
\begin{eqnarray*}
\lefteqn{\hspace{-1cm}P\left(Y_{T}>u_n \geq Y_{T+1},\ M^{(\bf{X})}_{T+2,\left[\frac{n}{k_nT}\right]T}>u_n\right)}\\[0.3cm]
& \leq & P\left(X_{T}>u_n \geq X_{T+1},\ M^{(\bf{X})}_{T+2,\left[\frac{n}{k_nT}\right]T}>u_n\right)\\[0.3cm]
&& + P\left( X_{T}>u_n \geq X_{T} \right)\\[0.3cm]
& \leq&   P\left( \bigvee_{j=0}^{T}  X_j>u_n \geq  X_{T+1},\ M^{(\bf{X})}_{T+2,\left[\frac{n}{k_nT}\right]T+T}>u_n \right).
\end{eqnarray*}
Due to (\ref{3.7}) with these last two upper bounds, we can now write
\begin{eqnarray}\label{3.8}
\lefteqn{n\frac{1}{T}\sum_{i=1}^{T}P\left(Y_i>u_n\geq M^{(\bf{Y})}_{i+1,i+T},\ M^{(\bf{Y})}_{i+T+1,\left[\frac{n}{k_nT}\right]T}>u_n\right)}\nonumber\\[0.3cm]
&\leq& \frac{2T-1}{T}\ nP\left(\bigvee_{j=0}^T X_j>u_n\geq X_{T+1},\ M^{(\bf{X})}_{T+2,\left[\frac{n}{k_nT}\right]T+T}>u_n\right),
\end{eqnarray}
Thus, condition $D_T^{(T+1)}(u_n)$ holds for $\{Y_n\}_{n\geq 1}$ since (\ref{3.8}) goes to zero as $n\to +\infty,$ form (\ref{3.6}).\edem\vspace{0.3cm}

Note that condition (\ref{3.6}) holds for $\{X_n\}_{n\geq 0}$ if it is a simple moving maxima sequence as defined in Example \ref{ex1}. In fact, $\{X_n\}_{n\geq 0}$ satisfies condition $D(u_n),$ for any sequence of real numbers $\{u_n\}_{n\geq 1},$ with $\alpha_{n,\ell_n}=0$ for $\ell_n\geq 3,$ since it is 2-dependent, and  if we consider  $r'_n=\left[\frac{n}{2k_n}\right],$ $u_n=nx,\, x>0,$ with $\tau=1/x,$ it holds
\begin{eqnarray*}
\lefteqn{\hspace{-1cm}nP\left(M_{0,2}^{(\bf{X})}>u_n\geq X_3, M^{(\bf{X})}_{4,\left[\frac{n}{2k_n}\right]2+2}>u_n\right)}\\[0.1cm]
&\leq& n\sum_{i=0}^2\sum_{j=4}^{2r'_n+2}P(X_i>u_n,X_j>u_n)\\[0.1cm]
&\leq& n(2r'_n+2)(1-F^2_Z(2u_n))\sum_{i=0}^2P(X_i>u_n)\\[0.1cm]
&\leq& 3n\left(\frac{n}{k_n}+2\right) (1-F^2_Z(2u_n))^2 \convn 0.
\end{eqnarray*}

We may then conclude that $\{Y_n\}_{n\geq 1},$ given in (\ref{1.2}) of Example \ref{ex1}, satisfies condition $D_2^{(3)}(u_n).$\vspace{0.3cm}

\begin{remark}
If the underlying stationary sequence $\{X_n\}_{n\geq 0}$ is an $m-$dependent sequences, {\it{i.e.}}, \linebreak $(X_{i_1},\ldots, X_{i_s})$ is independent of $(X_{j_1},\ldots,X_{j_t}),$ for all positive integers $s,t$ and $i_1<\ldots<i_s<j_1<\ldots<j_t$ with $j_1-i_s\geq m+1,$ then condition (\ref{3.6}) holds for some $T\leq m+1$ and $\{u_n\}_{n\geq 1}$ such that $n\overline{F}(u_n)\convn \tau_X.$
\end{remark}

\begin{remark} Condition  (\ref{3.6}) is implied by condition $$n\sum_{j=0}^T\sum_{k=T+2}^{\left[\frac{n}{k_nT}\right]T+T}P(X_j>u_n,X_{k-1}\leq u_n<X_k)\convn 0.$$
\end{remark}\vspace{0.3cm}

Under condition (\ref{3.6}) for $\{X_n\}_{n\geq 0}$ we shall see that it is possible to obtain  a relation between the extremal index of $\{Y_n\}_{n\geq 1},$ $\theta_{{\bf Y}},$ and the extremal index of $\{X_n\}_{n\geq 0},$ $\theta_{{\bf X}}.$
Before  we present a required lemma.

\begin{lemma}\label{lem1}$\quad$\newline\vspace{-0.5cm}
\begin{enumerate}
\item[{\bf{a)}}] For $i\in\{1,\ldots,T-1\}$ it holds
\begin{eqnarray} \label{3.9}
\lefteqn{\vspace{-1cm}P\left(Y_i>u_n, M^{(\bf{Y})}_{i+1,i+T} \leq u_n\right)}\\[0.3cm]
&=&P\left(X_i>u_n,M^{(\bf{X})}(S_i) \leq u_n,  M^{(\bf{Y})}(\overline{S}_i)\leq u_n \right)pP(A_i)\nonumber\\[0.3cm]
&&+P\left(\bigvee_{j=0}^{i-1}U_jX_j>u_n,  M^{(\bf{X})}(S_i) \leq u_n, M^{(\bf{Y})}(\overline{S}_i)\leq u_n \right)(1-p)P(A_i)\nonumber
\end{eqnarray}
\noindent with $S_i=\{i+1,\ldots,T\},$ $\overline{S}_i=\{T+1,\ldots,T+i\}$ and $\displaystyle{A_i:=\bigcap_{j\in S_i}\{U_j=1\}}.$\vspace{0.3cm}
\item[{\bf{b)}}] $P\left(Y_T>u_n, M^{(\bf{Y})}_{T+1,2T} \leq u_n\right)=P\left(X_T>u_n,  M^{(\bf{X})}_{T+1,2T} \leq u_n \right)p^{T-1}.$
\end{enumerate}
\end{lemma}
\bdem {\bf a)} We trivially have, for all $i\in\{1,\ldots,T-1\},$
\begin{eqnarray}
\lefteqn{\vspace{-1cm}P\left(Y_i>u_n, M^{(\bf{Y})}_{i+1,i+T} \leq u_n\right)}\nonumber\\[0.3cm]
&=&P\left(X_i>u_n,M^{(\bf{Y})}(S_i) \leq u_n,  M^{(\bf{Y})}(\overline{S}_i)\leq u_n \right)p\nonumber\\[0.3cm]
&&+P\left(\bigvee_{j=0}^{i-1}U_jX_j>u_n,  M^{(\bf{Y})}(S_i) \leq u_n, M^{(\bf{Y})}(\overline{S}_i)\leq u_n \right)(1-p).\nonumber
\end{eqnarray}
Consider now that $i\in\{1,\ldots,T-2\}.$  In order to deal with $M^{(\bf{Y})}(S_i)$, observe that  if there is some $m\in \{i+1,\ldots,T-1\}$ such that $U_m=0$, then  the probability (\ref{3.9}) involves the following probabilities
\begin{eqnarray}
& & P\left(X_i>u_n, \ldots,\bigvee_{j=0}^{m-1}U_jX_j\leq u_n, \ldots \Biggl|\  U_i=1,\ldots,U_m=0,\ldots \right)\nonumber \\[0.3cm]
&& \leq P(X_i <u_n \leq X_i)=0 \nonumber
\end{eqnarray}
and
\begin{eqnarray}
&&P\left(\bigvee_{j=0}^{i-1}U_jX_j>u_n, \ldots,\bigvee_{j=0}^{m-1}U_jX_j\leq u_n, \ldots \Biggl|\  U_i=0,\ldots,U_m=0,\ldots \right)\nonumber\\
&& \leq P\left(\bigvee_{j=0}^{i-1}U_jX_j < u_n \leq \bigvee_{j=0}^{i-1}U_jX_j\right)=0\nonumber
\end{eqnarray}
 since $m\geq i+1$. The proof follows immediately bacause $M^{(\bf{Y})}(S_i)$ leads directly to  $M^{(\bf{X})}(S_i)$ under the occurrence of $A_i$.
 For $i=T-1$ the result follows as well because $S_i=\{T\}$ and $M^{(\bf{Y})}(S_i)=X_T.$\vspace{0.5cm}

\noindent {\bf b)} Straightforward from the previously used arguments and the fact that $U_{kT}=1,$ $k\geq 0,$ $T\geq 1,$ almost surely.\edem\vspace{0.3cm}

\begin{theorem}\label{teo4}
If $\{X_n\}_{n\geq 0}$ satisfies condition (\ref{3.6}) for some $T\geq 2$ and $\{u_n\}_{n\geq 1}$  satisfying (\ref{3.4}), then
\begin{equation}\label{3.10}
\theta_{{\bf Y}}=\frac{\tau_{{\bf X}}}{\tau T}\theta_{{\bf X}}\left((T-1)p^{T}+p^{T-1}\right)+\displaystyle{\frac{1}{\tau T}\sum_{i=1}^{T-1} P_{i,T}},
\end{equation}
with
\begin{eqnarray*}
 P_{i,T} & = &\displaystyle{ \lim_{n\to
+\infty} }n \displaystyle{\sum_{ I_i \subsetneq \overline{S}_i}}P\left(X_{i}>u_n,\,  M^{(\bf{X})}(S_i \cup I_i) \leq u_n \right)p^{|I_i|+T-i}(1-p)^{i-|I_i|}\\[0.3cm]
 && +\displaystyle{ \lim_{n\to
+\infty} }n \displaystyle{\sum_{ J_i \subset \underline{S}_i} \sum_{I_i \subset \overline{S}_i}}P\left( M^{(\bf{X})}(J_i\cup \{0\})>u_n,\, M^{(\bf{X})}(S_i \cup I_i) \leq u_n \right)\times\\[0.3cm] &&\hspace{4cm} \times p^{|I_i|+|J_i|+T-i-1}(1-p)^{2i-|I_i|-|J_i|}
\end{eqnarray*}
where $J_1=\emptyset$,   $ S_i=\{i+1,\ldots,T\}, \, \overline{S}_i=\{T+1,\ldots,T+i\}$ for $i\in \{1,\ldots,T-1\}$ and, for $T \geq 3$ and  $i\in \{2,\ldots,T-1\}$,\, $ \underline{S}_i=\{1,\ldots,i-1\}$.
 \end{theorem}
\bdem Since $\{Y_n\}_{n\geq 1}$ satisfies $D_T^{(T+1)}(u_n)$, we deduce $\theta_{{\bf Y}}$ from (\ref{3.3}).
For $i=T$ we have, from Lemma \ref{lem1} b),
$$ P\left(Y_i>u_n,  M^{(\bf{Y})}_{i+1,i+T} \leq u_n \right)=P\left(X_{T}>u_n,  M^{(\bf{X})}_{T+1,2T} \leq u_n \right)p^{T-1}.$$
For any $i \in \{1,\ldots,T-1\},$ taking into account Lemma \ref{lem1} a), where $P(A_i)=p^{T-i-1},$ we can write
\begin{eqnarray*}
\lefteqn{P\left(Y_i>u_n, M^{(\bf{Y})}_{i+1,i+T} \leq u_n \right)}\nonumber\\[0.3cm]
& =&\displaystyle{\sum_{I_i \subset \overline{S}_i }}P\left(X_{i}>u_n,\, M^{(\bf{X})}(S_i) \leq u_n, M^{(\bf{X})}(I_i) \leq u_n \right)p^{|I_i|+T-i}(1-p)^{i-|I_i|}\nonumber\\[0.3cm]
&&  + \displaystyle{\sum_{ J_i \subset \underline{S}_i} \sum_{I_i \subset \overline{S}_i}}P\left(\max\{X_0, M^{(\bf{X})}( J_i )\}>u_n,\, M^{(\bf{X})}(S_i) \leq u_n, M^{(\bf{X})}(I_i) \leq u_n \right)\nonumber\\[0.3cm]
&& \hspace{4cm} \times p^{|I_i|+|J_i|+T-(i+1)}(1-p)^{i-|I_i|+1+(i-1)-|J_i|}\nonumber
\end{eqnarray*}
\begin{eqnarray*}
& = &P\left(X_{i}>u_n,M^{(\bf{X})}_{i+1,i+T} \leq u_n \right)p^{T} \\[0.3cm]
&& + \displaystyle{ \sum_{I_i \subsetneq \overline{S}_i}}P\left(X_{i}>u_n,\,  M^{(\bf{X})}(S_i \cup I_i) \leq u_n \right)p^{|I_i|+T-i}(1-p)^{i-|I_i|}\\[0.3cm]
&& + \displaystyle{\sum_{ J_i \subset \underline{S}_i}\sum_{I_i \subset \overline{S}_i}}P\left( M^{(\bf{X})}( J_i \cup \{0\} )>u_n,\, M^{(\bf{X})}(S_i \cup I_i) \leq u_n \right)\\[0.3cm]
&& \qquad \qquad \qquad \qquad \qquad \qquad \qquad \qquad \times p^{|I_i|+|J_i|+T-i-1}(1-p)^{2i-|I_i|-|J_i|}
\end{eqnarray*}

Now, due to the fact that condition (\ref{3.6}) implies condition $D^{(T+1)}(u_n)$ of  \cite{cher} and  $\{X_n\}_{n\geq 0}$ is a stationary sequence, the extremal index of $\{X_n\}_{n\geq 0}$ is given by $$\lim_{n\to
+\infty}P\left(M^{(\bf{X})}_{i+1,i+T} \leq u_n  \mid X_{i}>u_n \right) = \theta_{\bf{X}}.$$  Hence
$$\lim_{n\to
+\infty}nP\left(X_{i}>u_n,\,M^{(\bf{X})}_{i+1,i+T} \leq u_n  \right) =  \tau_{\bf{X}} \theta_{\bf{X}},$$
  which concludes the proof.
\edem \vspace{0.3cm}

We shall now apply the previous result in the computation of the extremal index of $\{Y_n\}_{n\geq 1},$ given in (\ref{1.1}), in two different scenarios. First, we shall consider a periodic control at instants multiple of $T=2,$ {\it{i.e.}} $\{Y_n\}_{n\geq 1}$ 2-periodic, and that the underlying stationary sequences $\{X_n\}_{n\geq 0}$ is the moving maxima of Example \ref{ex1}. Second, we consider a 3-periodic sequence $\{Y_n\}_{n\geq 1}$ with an underlying sequence  $\{X_n\}_{n\geq 0}$ a max-autoregressive sequence (ARMAX) of order one.

\begin{example}[{\bf{$\{X_n\}_{n \geq 0}$ moving maxima}}]$\quad$\newline
{\rm{As previously noted, sequence $\{Y_n\}_{n\geq 1}$ of Example \ref{ex1} satisfies condition $D_2^{(3)}(u_n)$ with $u_n=nx,\, x>0,$ therefore its extremal index, $\theta_{\bf Y},$ can be computed from (\ref{3.3}). This yields $\theta_{\bf Y}=\frac{1}{2}$ which is equal to $\theta_{\bf X}$ as expected.

Since $\{X_n\}_{n\geq 0}$ satisfies condition  (\ref{3.6}) for $u_n=nx,$ $x>0,$ the extremal index $\theta_{\bf Y}$  can also be calculated from (\ref{3.10}) of Theorem \ref{teo4}.  Note that in this case $\tau_{\bf X}=\tau$. Consequently
\begin{equation}\label{3.11}
\theta_{\bf Y}=\frac{1}{4}(p^2+p)+\frac{1}{2\tau}P_{1,2}
 \end{equation}
 with \begin{eqnarray}\label{3.11A}
P_{1,2}&=& \lim_{n\to +\infty} n \left(P(X_1>u_n,X_2\leq u_n) p(1-p)+\right.\\[0.3cm]
&&\left.+\ P(X_0>u_n,X_2\leq u_n)(1-p)^2+P(X_0>u_n,M_{2,3}^{({\bf{X}})}\leq u_n) p(1-p)\right).\nonumber
\end{eqnarray}
For the probabilities in $P_{1,2}$ we have
\begin{eqnarray*}
nP(X_1>u_n,X_2\leq u_n)&=&n(1-F_Z(2u_n)) F_Z^2(2u_n)\convn \frac{\tau}{2},\\[0.3cm]
nP(X_0>u_n,X_2\leq u_n)&=& n(1-F_Z^2(2u_n)) F_Z^2(2u_n)\convn \tau,\\[0.3cm]
nP(X_0>u_n,M_{2,3}^{({\bf{X}})}\leq u_n)&=& n(1-F_Z^2(2u_n)) F_Z^3(2u_n)\convn \tau.
\end{eqnarray*}
Taking this all into account, we finally obtain from (\ref{3.11}) that  $\theta_{\bf Y}=\frac{1}{2}.$}}
\end{example}\vspace{0.3cm}

\begin{example}[{\bf{$\{X_n\}_{n\geq 0}$ ARMAX of order one}}]$\quad$\newline

{\rm{If we have a solar thermal energy storage system where the temperature level in a tank is periodically controlled and eventually, for some reason, temperatures at certain time points are not retained, our model (\ref{1.1}) can be used to describe the temperature in such a situation. According to  Haslett (1979), the model defined by $$X_j =\beta X_{j-1} \vee (\alpha \beta X_{j-1} + Y_j), \quad j\geq 1, \quad 0\leq \alpha \leq 1, \quad 0<\beta<1,$$
can be used to describe the temperature level in a tank. The extremal behaviour of this first order ARMAX  storage model, for the  particular case $\alpha=0,$ was studied by Alpuim (1989) and in its multivariate version by  Ferreira and Ferreira (2013).

We shall now consider in (\ref{1.1}) $\{Y_n\}_{n\geq 1}$ to be a 3-periodic positive sequence and the underlying sequence $\{X_n\}_{n\geq 0}$ the first order ARMAX process of  Alpuim (1989),
$$ X_n=t\max\{X_{n-1},W_n\},\quad n\geq 1,$$ where $t\in]0,1[$ is a constant, $X_0$ is a positive random variable with d.f. $H_0,$ independent of the sequence of i.i.d. positive random variables $\{W_n\}_{n\geq 1}$ with d.f. $L.$

Let us assume that the Markovian sequence $\{X_n\}_{n\geq 1}$ is  stationary, {\it{i.e.}}, there exists $x>0$ such that $L(x/t)>0$ and $$0<\sum_{s=1}^{+\infty}(1-L(x/t^s))<+\infty,$$ as proved in  Alpuim (1989). Therefore, the non-degenerate d.f. $H$ of $X_n,$ $n\geq 0,$ satisfies the following equation $$L(x)=\frac{H(tx)}{H(x)},\quad x\geq \alpha(H)/t,$$ where $\alpha(H)=\inf\{x:H(x)>0\}\geq 0.$

It can be easily verified that, for $n\geq 1,$ $$X_n=\max\left\{t^nX_0,\max_{1\leq n\leq i}t^{i-n+1}Y_n,\right\}.$$
\vspace{0.3cm}

Sequence $\{X_n\}_{n\geq 1}$ satisfies the condition $D(u_n),$ for any sequence $\{u_n\}_{n\geq 1},$  because it is strong-mixing (see Alpuim 1988).

 If $H$ belongs to the max-domain of attraction of the Fr\'echet d.f with parameter $\alpha>0,$ then the normalized levels $\{u_n\}_{n\geq 1}$ for $\{X_n\}_{n\geq 1},$ ${\it{i.e.}},$ such that $n(1-H(u_n))\convn\tau_{{\bf{X}}}\geq 0,$  satisfy  $$n(1-H\left(u_n/t)\right)\convn\tau_{{\bf{X}}} t^{\alpha}.$$ In this case $\{X_n\}_{n\geq 1}$ has extremal index $\theta_{{\bf{X}}}=1-t^{\alpha}$ and
 \begin{eqnarray*}
  n\left(1-L(u_n/t)\right)&=& n\left(\frac{H(u_n/t)-H(u_n)}{H(u_n/t)}\right)\\[0.3cm]
  & = & \left(n(1-H(u_n)\right)-n\left(1-H(u_n/t))\right)
 \frac{1}{H(u_n/t)}\\[0.3cm]
  &\convn& \tau_{{\bf{X}}}(1-t^{\alpha})=\tau_{{\bf{X}}}\theta_{{\bf{X}}},
  \end{eqnarray*}
 (see Ferreira and Ferreira (2013) for further details). Similarly, we establish that
 \begin{eqnarray*}
n\left(H(u_n/t^{j+1})-H(u_n)\right)& \convn& \tau_{{\bf{X}}}\theta_{{\bf{X}}}(1+t^{\alpha}+\ldots+t^{j\alpha}),\, j \geq 0,\\[0.3cm]
n\left(1-L\left(u_n/t^{j}\right)\right) &\convn& t^{(j-1) \alpha} \tau_{{\bf{X}}}\theta_{{\bf{X}}}, \, j \geq 1,\\[0.3cm]
n\left(1-L^j\left(u_n/t\right)L\left(u_n/t^{2}\right)\right) &\convn&  \tau_{{\bf{X}}}\theta_{{\bf{X}}}(j+t^{\alpha}),\, j \geq 1.
	\end{eqnarray*}

 Condition  (\ref{3.6}) holds for  $\{X_n\}_{n\geq 0}$, for the normalized levels $\{u_n\}_{n\geq 1}$ and any $T\geq 2,$ since\vspace{0.2cm}
\begin{eqnarray*}
\lefteqn{\hspace{-2cm}nP\left(\bigvee_{j=0}^T X_j>u_n\geq X_{T+1},\ M^{(\bf{X})}_{T+2,r'nT+T}>u_n\right)}\\[0.3cm]
&\leq& n\frac{1}{T}\sum_{j=0}^{T}\sum_{s=T+1}^{r'_nT+T} P(X_j>u_n,X_s \leq u_n <X_{s+1})%\\[0.3cm]
\end{eqnarray*}
\begin{eqnarray*}
&\leq&n\frac{1}{T}\sum_{j=0}^{T}\sum_{s=T+1}^{r'_nT+T} P(X_j>u_n,X_s \leq u_n, W_{s+1}>u_n/t)\\[0.3cm]
&\leq& n\frac{1}{T}\sum_{j=0}^{T}\sum_{s=T+1}^{r'_nT+T} P(X_j>u_n)P( W_{s+1}>u_n/t)\\[0.3cm]
&\leq &
n\frac{T+1}{T}\left(\frac{n}{k_n}+T\right)(1-H(u_n)) \left(1-L(u_n/t)\right)\\[0.3cm]
&  \convn& 0,
\end{eqnarray*}
 for any positive integer sequence $\{k_n\}$ such that $k_n\convn + \infty$ and\linebreak  $k_n/n \convn 0.$

The validation of condition (\ref{3.6}) for  $\{X_n\}_{n\geq 0}$ guarantees that condition $D^{(4)}_3(u_n)$ holds for $\{Y_n\}_{n\geq 1}$ (Theorem \ref{teo3}) and therefore its extremal index can be computed from the expression given in Theorem \ref{teo4}.
Indeed, since all the factors with products containing  $L^j\left(u_n / t^m\right)$ tend to 1, for all $j,\,m \geq 1$, we have
\begin{eqnarray*}
  \displaystyle{\sum_{i=1}^{T-1} P_{i,T}}& = & \lim_{n\to
+\infty} \Big\{ n\left(H(u_n/t) - H(u_n)\right)(p+p^2-2p^3) +  \\[0.3cm]
&&+\ n\left( H(u_n/t^2)-H(u_n)\right)p(1-p) +\\[0.3cm]
&&+\ n\left( H(u_n/t^3) -H(u_n)\right)(1-p)^2 +  \\[0.3cm]
&&+\ n\left( 1-L^2(u_n/t)L(u_n/t^2)\right)p(1-p)\Big\}\\[0.3cm]
& = & \tau_{{\bf{X}}}\theta_{{\bf{X}}} \Big\{ \left(p+p^2-2p^3\right)+(2-\theta_{{\bf{X}}})p(1-p)+(3-\theta_{{\bf{X}}})p(1-p)\\[0.3cm]
 &  &  +\ (3-3\theta_{{\bf{X}}}+\theta^2_{{\bf{X}}})\left(1-p\right)^2 \Big\}.\\
\end{eqnarray*}

Furthermore,  $G_1(x)=H(x)$ and $G_2(x)=(1-p)H(x)-pP(X_1 \leq x,\, X_0 \leq x),$ leading (\ref{3.4})  to
$\tau_1=\tau_{\bf{X}}$ and $\tau_2=\tau_\textbf{X}+p \theta_\textbf{X} \tau_\textbf{X}$. Hence, by  (\ref{3.5}), it holds $$\tau=\tau_\textbf{X}\left(1+\frac{1}{3}p(1-p) \theta_\textbf{X}\right).$$  As a consequence, (\ref{3.10}) becomes
\begin{eqnarray*}
\theta_{{\bf Y}}& =& \frac{\tau_{{\bf X}}}{\tau T}\theta_{{\bf X}}\left(2p^3+p^2\right)+\frac{\tau_{{\bf X}}}{\tau T}\theta_{{\bf X}} \left(3-2p^3-p^2+\theta_{{\bf X}}(-3+4p-p^2)\right. +\\
&&+\ \left. \theta^2_{{\bf X}}(1-p)^2\right)\\[0.3cm]
& =& \frac{3\theta_\textbf{X}+ \theta_\textbf{X}^2(-3+4p-p^2)+ \theta_\textbf{X}^3(1-p)^2}{3+p(1-p) \theta_\textbf{X}}.
\end{eqnarray*}

Figure \ref{fig2} shows the effect that the extremal index of the underlying sequence $\theta_\textbf{X}\in ]0,1]$ and the parameter $p\in]0,1[$ have on the extremal index of $\{Y_n\}_{n\geq 1},$ $\theta_{{\bf Y}}.$ As we can see, when $p$ is very close to one we get $\theta_{{\bf Y}}\simeq \theta_\textbf{X},$ since in this case there are almost no missing values and so $Y_n=X_n,$ $n\geq 1,$ almost surely. When $p$ is very close to zero, almost all values of the underlying sequence are missing, except for the values at instants multiple of three, since $T=3.$ In this case, sequence  $\{Y_n\}_{n\geq 1}$ will have the following form $X_0,X_0,X_3,X_3,X_3,X_6,X_6,X_6,..., X_{sT},X_{sT},X_{sT}\ldots,$ and so if $\theta_\textbf{X}\simeq 1$ ($\theta_\textbf{X}=1$ occurs, for example, for i.i.d  sequences), exceedances of high levels will form clusters of mean size approximately $T=3,$ yielding an extremal index approximately equal to $1/3$ as observed in Figure \ref{fig2}.
\begin{figure}[!ht]\label{fig2}
 \begin{center}
 \includegraphics[scale=0.45]{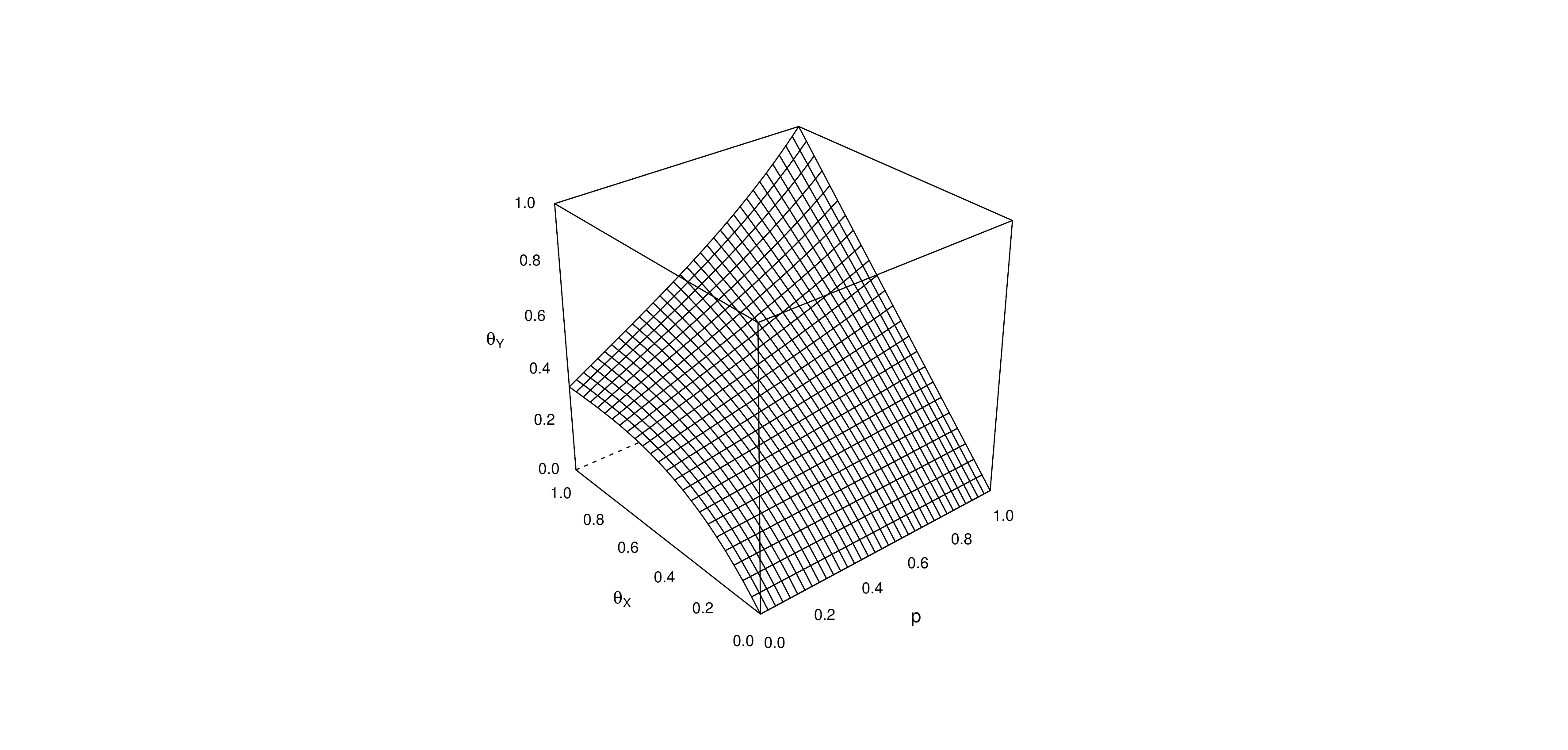}
 \caption{$\theta_Y$ as a function of $\theta_X$ and $p$ for an ARMAX underlying sequence}
 \end{center}
 \end{figure}

Moreover, considering that the periodic control took place at instants multiple of two, $T=2,$ then, with $P_{1,2}$ given by (\ref{3.11A}) we would obtain, from Theorem \ref{teo4},
\begin{eqnarray*}
\theta_{{\bf Y}}& =& \frac{\theta_\textbf{X}}{2}(p^2+p)+ \frac{\theta_\textbf{X}}{2}\left((2-\theta_{{\bf X}})(1-p)+p(1-p)\right)\\[0.3cm]
& =& \theta_\textbf{X}+ \theta_\textbf{X}^2\frac{p-1}{2},
\end{eqnarray*}
since it was proved that condition (\ref{3.6}) also holds for the underlying ARMAX sequence of order one when $T=2.$}}
  \end{example}

Imposing a stronger condition on the behaviour of the underlying stationary sequence $\{X_n\}_{n\geq 0}$ than condition  (\ref{3.6}), we obtain, with similar arguments as used to prove Theorem \ref{teo3},  the validation of condition  $D_T^{(T)}(u_n)$ for  $\{Y_n\}_{n\geq 1}$ and consequently also of $D_T^{(T+1)}(u_n),$ as stated in the next result.\vspace{0.3cm}

\begin{theorem}\label{teo5}
If for any sequence $\{u_n\}_{n\geq 1},$ the sequence $\{X_n\}_{n\geq 0}$ satisfies
\begin{equation}\label{3.12}
nP\left(\bigvee_{j=0}^T X_j>u_n,\ M^{(\bf{X})}_{T+1,\left[\frac{n}{k_nT}\right]T+T}>u_n\right)\convn 0
\end{equation}
 then condition $D_T^{(T)}(u_n)$ holds for the sequence $\{Y_n\}_{n\geq 1}.$
\end{theorem}\vspace{0.3cm}

Condition (\ref{3.12}) is indeed more demanding than condition (\ref{3.6}) as we can verify  with the moving maxima sequence $\{X_n\}_{n\geq 0}$  defined in Example \ref{ex1}. In this case condition (\ref{3.12}) with $T=2$ does not hold  since
\begin{eqnarray*}
\lefteqn{\hspace{-1cm}nP\left(M_{0,2}^{(\bf{X})}>u_n, M^{(\bf{X})}_{3,\left[\frac{n}{2k_n}\right]2+2}>u_n\right)}\\[0.3cm]
&>& nP(X_2>u_n,X_3>u_n)\\[0.3cm]
&= & n(1-F_Z(2u_n))\left(1+F_Z(2u_n)-F_Z^2(2u_n)\right) \convn \tau >0.
\end{eqnarray*}
Furthermore, the associated sequence $\{Y_n\}_{n\geq 1},$  also does not satisfy condition $D^{(2)}_2(u_n)$ since
\begin{eqnarray*}
\lefteqn{\hspace{-1cm}nP\left(Y_1>u_n\geq Y_2, M^{({\bf{Y}})}_{3,\left[\frac{n}{2k_n}\right]2}>u_n\right)}\\[0.3cm]
&\geq& nP(Y_1>u_n\geq Y_2, Y_3>u_n)\\[0.3cm]
&\geq & nP(X_0>u_n\geq X_2, X_0>u_n) \convn \tau >0.
\end{eqnarray*}\vspace{0.3cm}

\begin{remark}
Condition (\ref{3.12}) is implied by the following condition $$n\sum_{j=0}^T\sum_{k=T+1}^{\left[\frac{n}{k_nT}\right]T+T}P(X_j>u_n, X_k>u_n)\convn 0.$$
\end{remark}\vspace{0.3cm}

We now present an analogous result to Theorem \ref{teo4}, relating $\theta_{\bf{Y}}$ and $\theta_{\bf{X}}$ under the hypothesis that $\{X_n\}_{n\geq 0}$ satisfies condition (\ref{3.12}) and consequently $\{Y_n\}_{n\geq 1}$ satisfies condition $D_T^{(T)}(u_n).$  The main difference can be found in the first term of the expression for $\theta_{\bf{Y}}$ and in the definition of $\overline{S}_i.$ \vspace{0.3cm}

\begin{theorem}\label{teo6}
If $\{X_n\}_{n\geq 0}$ satisfies condition (\ref{3.12}), for some $T\geq 2$ and $\{u_n\}_{n\geq 1}$  satisfying (\ref{3.4}), then $\{Y_n\}_{n\geq 1}$ satisfies condition $D_T^{(T)}(u_n)$  and  $$\theta_{{\bf Y}}=\frac{\tau_{{\bf X}}}{\tau }\theta_{{\bf X}}p^{T-1}+\displaystyle{\frac{1}{\tau T}\sum_{i=1}^{T-1} P^*_{i,T}}  $$
with
$$\begin{array}{lll}
 P^*_{i,T} & = & \displaystyle{\lim_{n\to +\infty}}n\left(\displaystyle{\sum_{I_i \subsetneq \overline{S}_i }}P\left(X_{i}>u_n,\,  M^{(\bf{X})}(S_i \cup I_i) \leq u_n \right)\right.\times\\ && \qquad \qquad \times\ p^{|I_i|+T-i}(1-p)^{i-1-|I_i|}+\\[0.3cm]
 &  &+\ \left.\displaystyle{\sum_{J_i \subset \underline{S}_i}\sum_{I_i \subset \overline{S}_i }}P\left( M^{(\bf{X})}(J_i\cup \{0\})>u_n,\, M^{(\bf{X})}(S_i \cup I_i) \leq u_n \right)\right.\times\\
 && \qquad \qquad \times \ \left. p^{|I_i|+|J_i|+T-i-1}(1-p)^{2i-1-|I_i|-|J_i|}\right),
\end{array}$$ where
$I_1=J_1= \emptyset,\,$   $S_i=\{i+1,\ldots,T\}$ for any $i \,\in \, \{1,...,T-1\}$ and, for $T\geq 3$ and any $i \,\in \, \{2,...,T-1\}$ $ \overline{S}_i = \{T+1,\ldots,T+i-1\}$  and $ \underline{S}_i=\{1,...,i-1\}$.
\end{theorem}\vspace{0.3cm}

We observe that any i.i.d. positive sequence $\{ X_n\}$
 trivially satisfies  (\ref{3.12}) for normalized levels $u_n$, and therefore sequence $\{Y_n\}$ satisfies condition $D_T^{(T)}(u_n)$.
Then, considering for instance $T=3$  and proceeding as before,  we get $\tau_2= \tau_\textbf{X}(1+p)$, $\,\, \tau=\frac{1}{3}\tau_\textbf{X}(3+p-p^2)$ and
  \begin{eqnarray*}
P^*_{i,2}& = &\lim_{n\to
+\infty} n (1-F(u_n))(1-p^2) + \lim_{n\to
+\infty} n (1-F^2(u_n))(p-p^2) \\[0.3cm]
& = & \tau_{\bf{X}}(1+2p-3p^2).
\end{eqnarray*}
 Hence
$$\theta_Y =\frac{\tau_{\bf{X}}}{\tau}p^2+\frac{1}{3 \tau}\tau_{\bf{X}}(1+2p-3p^2) =\frac{1+2p}{{3+p(1-p)}}.$$

%\begin{acknowledgements}
%\end{acknowledgements}

\end{document}